\documentclass[12pt]{article}
\usepackage{amsmath, amsthm, amssymb}
\usepackage[english]{babel}
\usepackage{graphicx}
\usepackage{float}
\usepackage{url}
\usepackage{caption}

\begin{document}

\title{\bf How do curved spheres intersect in 3-space?}
\author{Sergey Avvakumov, \texttt{s.avvakumov@gmail.com}}
\maketitle

\begin{abstract}
The following problem was proposed in 2010 by S. Lando.

Let $M$ and $N$ be two unions of the same number of disjoint circles in a sphere.
Do there always exist two spheres in 3-space such that their intersection is transversal and is a union
of disjoint circles that is situated as $M$ in one sphere and as $N$ in the other? Union $M'$ of disjoint circles is {\it situated} in one sphere as union $M$ of disjoint circles in the other sphere if there is a homeomorphism between these two spheres which maps $M'$ to $M$.

We prove (by giving an explicit example) that the answer to this problem is ``no''. We also prove a necessary and sufficient condition on $M$ and $N$ for existing of such intersecting spheres.
This result can be restated in terms of graphs. Such restatement allows for a trivial brute-force algorithm checking the condition for any given $M$ and $N$. It is an open question if a faster algorithm exist.
\end{abstract}

\newtheorem{thm1}{Theorem}
\newtheorem{thm2}[thm1]{Theorem}
\newtheorem*{thmExt}{Embedding Extension Theorem}
\newtheorem{clm1}{Claim}

\theoremstyle{definition}
\newtheorem{pbm1}{Open problem}
\newtheorem{pbm2}[pbm1]{Open problem}

\theoremstyle{remark}
\newtheorem{exm1}{Example}
\newtheorem{exm2}[exm1]{Example}
\newtheorem*{borr}{Borromean rings}
\newtheorem{cor2}{Corollary}

\subsection*{The Lando Problem}
We work entirely in the piecewise-linear (PL) category \footnote{A {\it PL circle} or {\it circle} is a closed broken line (polygon) without self-intersections in 3-space. A {\it PL sphere} or {\it sphere} is a polyhedron in 3-space (more precisely, 2-dimensional surface of the polyhedron), which is split
into several parts by any circle lying on the polyhedron, i.e. is a polyhedron homeomorphic to $S^2$.}.

Suppose $M$ and $M'$ are the unions of the same number of disjoint circles in spheres $S$ and $S'$. Then {\it $M$ is situated in $S$ as $M'$ in $S'$} if there is a homeomorphism $f:S\rightarrow S'$ such that $f(M)=M'$.

The following problem suggested by S. Lando was one of the (unsolved) problems at the Moscow State University mathematical tournament for students and young professors 2010 (\cite{Lando}, problem MB-8).

\smallskip
{\it Let $M$ and $N$ be two unions of the same number of disjoint circles in a sphere.
Do there exist two spheres in 3-space whose intersection is transversal and is a union
of disjoint circles that is situated as $M$ in one sphere and as $N$ in the other?}
\smallskip

This problem appeared in the discussion of related papers \cite{Rukhovich}, \cite{Hirasa}, \cite{Nowik}.

In this paper we prove that the answer to Lando problem is ``no'' by giving an explicit example.

	\begin{figure}[H]
	\centering
	\includegraphics[width=100mm]{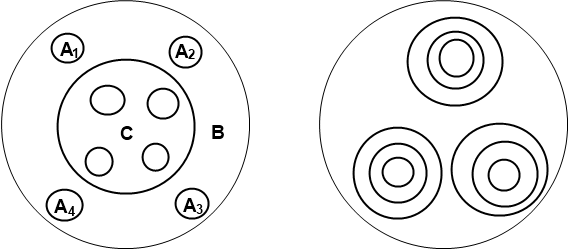}
	\caption{Two unions of $M$ (left) and $N$ (right) of $9$ circles.}
	\label{fig:CounterEx}
	\end{figure}
	
\begin{thm1}[an example]
\label{thm1}
Let $M$ and $N$ be two unions of $9$ disjoint circles in $S^2$ shown in Fig.~\ref{fig:CounterEx}. Then there are no two spheres in 3-space whose intersection is transversal and is a union of $9$ disjoint circles that is situated as $M$ in one sphere and as $N$ in the other. 
\end{thm1}

	\begin{figure}[H]
	\centering
	\includegraphics[width=100mm]{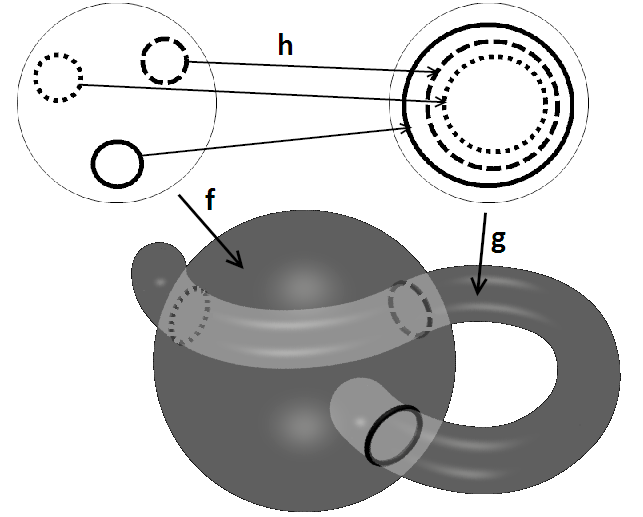}
	\caption{Bijection $h$ between two sets of three circles is realized by PL embeddings $f$, $g$.}
	\label{fig:BijectionEx}
	\end{figure}
	
In Theorem \ref{thm2} (see below) we describe all the collections of circles which can be realized by two intersecting spheres. The precise meaning of the word ``realized'' is defined in the following paragraph.

Assume that $M$ and $N$ are two unions of disjoint circles in sphere $S^2$. Suppose there exists PL embeddings\footnote{Map $f:A\rightarrow B$ is piecewise linear if $f$ is a simplicial map for {\it some} simplicial decompositions of $A$ and $B$} $f:S^2\hookrightarrow {\mathbb R^3}$ and $g:S^2\hookrightarrow {\mathbb R^3}$ such that intersection $f(S^2)\cap g(S^2)$ is transversal and $f(S^2)\cap g(S^2)=f(M)=g(N)$. These embeddings induce a bijection $h$ between sets of circles of $M$ and of $N$ (for circles $m\subset M$ and $n\subset N$ let $h(m)=n$ if $f(m)=g(n)$). Equivalently we may number circles of $M$ and of $N$ by $1,\dots,k$ so that two circles corresponding to the same circle of $f(S^2)\cap g(S^2)$ have the same number. We say that $f$, $g$ {\it realize} $h$ (Fig.~\ref{fig:BijectionEx}).

Theorem \ref{thm2} (see below) gives a necessary and sufficient condition for the realizability of a bijection. In particular Theorem \ref{thm2} can be used to prove Theorem \ref{thm1}. The following simple example shows that not every bijection is realizable. 

	\begin{figure}[H]
	\centering
	\includegraphics[width=100mm]{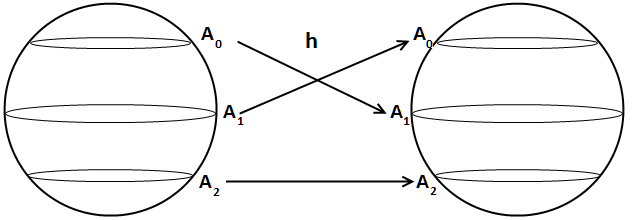}
	\caption{Circles $A_0,A_1,A_2$, bijection $h$.}
	\label{fig:NumberedEx}
	\end{figure}
	
\begin{exm1}
\label{exm1}
Let $A_0,A_1,A_2$ be the circles situated in $S^2$ as shown in the Fig.~\ref{fig:NumberedEx}. Let $M=N=A_0\cup A_1\cup A_2$. Let $h$ be a bijection defined by $h(A_0)=A_1$, $h(A_1)=A_0$ and $h(A_2)=A_2$. Then $h$ is not realizable.
\end{exm1}

The proof of Example \ref{exm1} (see ``Proofs'') demonstrates some of the ideas used in the proof of Theorem \ref{thm2}.

Let us introduce definitions necessary to state Theorem \ref{thm2}.

Let $M$ and $N$ be two unions (not necessary nonempty) of disjoint circles in sphere $S^2$.
Color connected components of $S^2-N$ in black and white so that adjacent components have different colors.
Union $M$ is {\it on one side} (in this sphere) of $N$ if $M$ is contained in the union of same colored components
of $S^2-N$. Unions $M$ and $N$ are {\it unlinked} (in this sphere) if $M$ is on one side of $N$ and $N$ is on one side of $M$. Equivalently unions $M$ and $N$ are {\it unlinked} (in $S^2$) if $[M]=0$ in $H_1(S^2-N;{\mathbb Z}_2)$ and $[N]=0$ in $H_1(S^2-M;{\mathbb Z}_2)$.

	\begin{figure}[H]
	\centering
	\includegraphics[width=120mm]{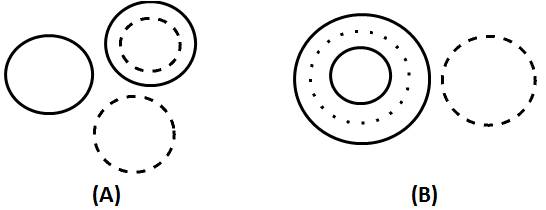}
	\caption{(A) $M$ (solid) is on one side of $N$ (dashed) while $N$ is not on one side of $M$.\\(B) $M$ (solid) and $N$ (dashed), $N$ and $P$ (dotted) are unlinked, but $M$ and $P$ are not unlinked.}
	\label{fig:UnlinkProperties}
	\end{figure}
	
Unions $M$ and $N$ are always unlinked if $M$ or $N$ is empty. If $M$ is on one side of $N$ then $N$ is not necessary on one side of $M$ (Fig.~\ref{fig:UnlinkProperties} left). Unlinkedness is not transitive. That is, if $M$ and $N$, $N$ and $P$ are unlinked, then $M$ and $P$ are not necessarily unlinked (Fig.~\ref{fig:UnlinkProperties} right).

Let $M$ be a union of disjoint circles in sphere $S$. Suppose $A$ is a connected component of $S-M$. Denote by $\partial A$ the boundary of the closure of $A$.

\begin{thm2}
\label{thm2}
Let $M$ and $N$ be two unions of the same number of disjoint circles in $S^2$. Let $h$ be a bijection between sets of circles of $M$ and of $N$. Color connected components of $S^2-M$ in two colors so that any two same colored components are not adjacent. Then $h$ is realizable if and only if $h(\partial A)$ and $h(\partial B)$ are unlinked for each two same-colored components $A$ and $B$ of $S^2-M$.
\end{thm2}

We say that sphere with holes $P$ is {\it properly embedded} in $D^3$ if $\partial P \subset \partial D^3$ and the interior of $P$ lies in the interior of $D^3$. Theorem \ref{thm2} is proved using the following:
\begin{thmExt}
\label{thmExt}
Let $M_1,\ldots,M_m$ be unions of disjoint circles in the sphere $S^2=\partial D^3$. Then there exist properly embedded in $D^3$ disjoint spheres with holes $P_1,\ldots,P_m$ such that $\partial P_i=M_i$ for each $i=1,\dots,m$ if and only if $M_1,\ldots,M_m$ are pairwise unlinked.
\end{thmExt}

Embedding Extension Theorem immediately implies the following:
\begin{cor2}
\label{cor2}
Let $M_1,\ldots,M_m$ be unions of disjoint circles in the sphere $S^2=\partial D^3$. Suppose that for every $i,j$ there exist properly embedded in $D^3$ disjoint spheres with holes $P'_i,P'_j$ such that $\partial P'_i=M_i, \partial P'_j=M_j$. Then there exist properly embedded in $D^3$ disjoint spheres with holes $P_1,\ldots,P_m$ such that $\partial P_i=M_i$ for each $i=1,\dots,m$.
\end{cor2}

\begin{figure}[H]
	\centering
	\includegraphics[width=90mm]{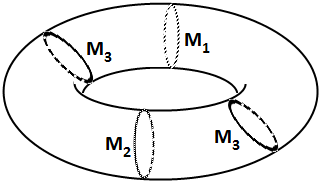}
	\caption{Unions $M_1$ and $M_2$ consists of one meridian each and $M_3$ consists of two meridians.}
	\label{fig:BorromeanRings}
	\end{figure}

Note that analogous statement is false for all closed orientable 2-surfaces other than $S^2$. For instance:
\begin{exm2}
\label{exm2}
Let $M_1,M_2,M_3$ be unions of disjoint circles in the standard torus $T^2\subset{\mathbb R^3}$. Let $M_1$ and $M_2$ be a single meridian each and let $M_3$ be a union of two meridians (Fig.~\ref{fig:BorromeanRings}). Then for every $i,j$ there exist disjoint spheres with holes $P'_i,P'_j$ whose interiors are inside $T^2$ and such that $\partial P'_i=M_i, \partial P'_j=M_j$. But there are no disjoint spheres with holes $P_1,P_2,P_3$ whose interiors are inside $T^2$ and such that $\partial P_1=M_1,\partial P_2=M_2,\partial P_3=M_3$.
\end{exm2}

This example is similar to the famous Borromean rings example stated in the following way:
\begin{borr}
\label{borr}
Let $S^1_1,S^1_2,S^1_3$ be the Borromean rings in $S^3=\partial D^4$. Then for every $i,j$ there exist properly embedded in $D^4$ disjoint disks $D'^2_i,D'^2_j$ such that $\partial D'^2_i=S^1_i, \partial D'^2_j=S^1_j$. But there are no properly embedded in $D^4$ disjoint disks $D^2_1,D^2_2,D^2_3$ such that $\partial D^2_1=S^1_1,\partial D^2_2=S^1_2,\partial D^2_3=S^1_3$.
\end{borr}

\subsection*{Relation to graphs}
Suppose that $M$ is a union of disjoint circles in sphere $S^2$.
Define (``dual to $M$'') graph $G=G(S^2,M)$ as follows.
The vertices are the connected components of $S^2-M$. Two vertices are connected
by an edge if the corresponding connected components are neighbors.

	\begin{figure}[H]
	\centering
	\includegraphics[width=120mm]{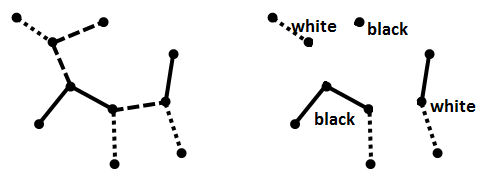}
	\caption{Graph $G$ (left), $p$ (solid edges), $q$ (dashed edges), complement in $G$ to the interiors of edges of $q$ (right). Set $p$ is not on one side of $q$ since edges of $p$ are both in black and white connected components. Thus $p$ and $q$ are not unlinked.}
	\label{fig:UnlinkEdgeEx}
	\end{figure}
	
The definition of unlinked unions of circles can also be restated in terms of graphs.
Let $p$ and $q$ be two sets of edges of a tree $G$.
Color connected components of the complement in $G$ to the interiors of edges of $q$ in black and white so that adjacent components have different colors.
The set $p$ is {\it on the same side} of $q$ (in this tree $G$)  if
$p$ is contained in the union of same-colored connected components of $G-q$
(or, equivalently, if $p\cap q=\emptyset$ and for each two vertices of edges of $p$ there is a path in the tree connecting these two points, and containing an even number of edges of $q$).
Sets $p$ and $q$ are {\it unlinked} (in this tree) if $p$ is on the same side of $q$ and $q$ is on the same side of $p$ (for example see Fig.~\ref{fig:UnlinkEdgeEx}).

Let $G$ and $H$ be two trees with the same number of edges. Color vertices of $G$ in two colors so that any two same colored vertices are not adjacent. Bijection $h$ between the sets of edges of $G$ and $H$ is called {\it realizing} if $h(\delta A)$ \footnote{$\delta A$ is a set of all edges incident to $A$} and $h(\delta B)$ are unlinked (in $H$) for each two same-colored vertices $A$ and $B$ of $G$.

Instead of a union of disjoint circles in a sphere let us consider its dual graph. Theorem \ref{thm2} implies that a bijection between two sets of circles is realizable if and only if the corresponding bijection between the sets of edges of dual graphs is realizing.

Let $G$ and $H$ be two trees with $k$ edges each. Given a bijection $h$ between the sets of edges of $G$ and $H$ we can check algorithmically in at most $O(k^2)$ time if $h$ is realizing. So, there is a brute-force algorithm which finds a realizing bijection (if any) in $O(k^2k!)$ time. We don't know if the more efficient algorithm exists. More precisely there is the following open problem:

\begin{pbm1}
\label{pbm1}
Is there a ``fast'' algorithm, which takes as input two arbitrary trees $G$ and $H$ with $k$ edges each and produces as output a realizing bijection (if any) between the sets of edges of $G$ and $H$?
\end{pbm1}

\begin{pbm2}
\label{pbm2}
Is there a tree $G$ such that there is no realizing bijection between the sets of edges of $G$ and $H$, where $H$ is the path graph with the same number of edges as $G$?
\end{pbm2}

\subsection*{Proofs}

\begin{proof}[Proof of the Example \ref{exm1}]
Assume to the contrary that there are PL embeddings $f:S^2\hookrightarrow {\mathbb R^3}$ and 
$g:S^2\hookrightarrow {\mathbb R^3}$ realizing $h$.

	\begin{figure}[H]
	\centering
	\includegraphics[width=130mm]{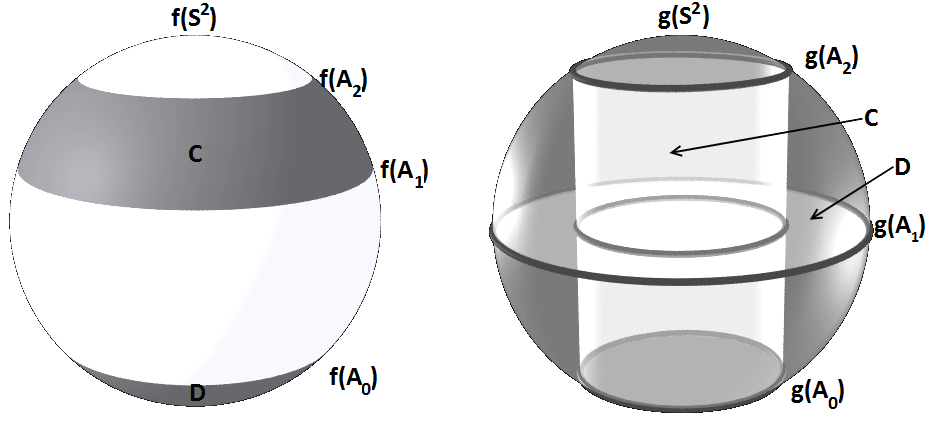}
	\caption{Circles $f(A_0)=g(A_1)$, $f(A_1)=g(A_0)$, $f(A_2)=g(A_2)$.}
	\label{fig:NumberedProof}
	\end{figure}
	
Denote by $D$ the disk in $f(S^2)-g(S^2)$ bounded by $f(A_0)$ (Fig.~\ref{fig:NumberedProof}). Denote by $C$ the cylinder in $f(S^2)-g(S^2)$ bounded by $f(A_1)$ and $f(A_2)$. Clearly $C$ and $D$ lie in 3-space on the same side of sphere $g(S^2)$. Circles $f(A_1)=g(h(A_1))=g(A_0)$ and $f(A_2)=g(h(A_2))=g(A_2)$ lie in sphere $g(S^2)$ on the different sides of the circle $f(A_0)=g(h(A_0))=g(A_1)$. So $C$ intersects $D$. This contradicts to the assumption that $f$ is an embedding.
\end{proof}

\begin{proof}[Proof of Theorem \ref{thm1}]
Assume to the contrary that there is a bijection $h$ between sets of circles of $M$ and of $N$ and PL embeddings $f:S^2\hookrightarrow {\mathbb R^3}$ and $g:S^2\hookrightarrow {\mathbb R^3}$ realizing $h$.

Denote the connected components of $f(S^2)-g(S^2)$ as shown in Fig.~\ref{fig:CounterEx} (left).

Consider disks $A_1,\ldots,A_4\subset f(S^2)$. Without loss of generality we may assume that their interiors lie inside $g(S^2)$. Then the interior of component $C\subset f(S^2)$ lies inside $g(S^2)$ as well (since the intersection $f(S^2)\cap g(S^2)$ is transversal). Since $C,A_1,\ldots,A_4$ are disjoint, $C$ lies completely in one of the connected components of ${\mathbb R^3}-g(S^2)\cup\bigsqcup A_i$. So all the $5$ circles of $\partial C$ lie in the same connected component of $g(S^2)-\bigsqcup\partial A_i$ (this argument is generalized in the proof of Claim \ref{clm1} below).

	\begin{figure}[H]
	\centering
	\includegraphics[width=50mm]{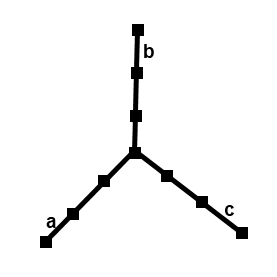}
	\caption{Graph $G(S^2,N)$.}
	\label{fig:TriodeTree}
	\end{figure}
	
Let us restate the previous statement in terms of graph $G(S^2,N)$ (Fig.~\ref{fig:TriodeTree}). Denote by $G(C)$ the union of $5$ edges of $G(S^2,N)$ corresponding to the circles of $\partial C$. Then $G(C)$ lies completely in one of the connected components of the compliment of $G(S^2,N)$ to the $4$ edges corresponding to the circles of $\bigsqcup\partial A_i$. Since $G(S^2,N)$ has only $9$ edges this means that $G(C)$ is a subtree of $G(S^2,N)$. Denote by $G(B)$ the union of $5$ edges of $G(S^2,N)$ corresponding to the circles of $\partial B$. Likewise, $G(B)$ is a subtree of $G(S^2,N)$.

Since $G(B)\cup G(C)=G(S^2,N)$, at least two of the three edges $a,b,c$ of $G(S^2,N)$ (Fig.~\ref{fig:TriodeTree}) belong to one of subtrees $G(B)$ or $G(C)$. Without loss of generality we may assume that $a,b\in G(B)$. But any subtree of $G(S^2,N)$ containing both $a$ and $b$ has at least $6$ edges while $G(B)$ has only $5$ edges. This contradicts the initial assumption.
\end{proof}

\begin{proof}[Proof of the ``only if'' part of Theorem \ref{thm2}]
Let $A$ and $B$ be two same colored components of $S^2-M$. Then $f(A)$ and $f(B)$ lie on the same side of $g(S^2)$. So by the ``only if'' part of Embedding Extension Theorem $\partial f(A)$ and $\partial f(B)$ are unlinked in $g(S^2)$. Then $g^{-1}(\partial f(A))=h(\partial A)$ and $g^{-1}(\partial f(B))=h(\partial B)$ are unlinked in $S^2$.
\end{proof}

\begin{proof}[Proof of the ``if'' part of Theorem \ref{thm2}]
Let $g:S^2\hookrightarrow {\mathbb R^3}$ be any PL embedding. We define a PL embedding $f:S^2\hookrightarrow {\mathbb R^3}$ such that $f$, $g$ realize $h$ by defining $f(A)$ for every connected component $A$ of $S^2-M$.

Color connected components of $S^2-M$ in black and white so that any two same colored components are not adjacent.

Let $P_1,\ldots,P_m$ be the white components of $S^2-M$. By the assumption of the Theorem $h(\partial P_1),\ldots,h(\partial P_m)$ are pairwise unlinked in $S^2$. 
So $g(h(\partial P_1)),\ldots,g(h(\partial P_m))$ are pairwise unlinked in $g(S^2)$. 
By the ``if'' part of Embedding Extension Theorem there exist disjoint spheres with holes $P'_1,\ldots,P'_m$
 whose interiors are inside $g(S^2)$ and such that $\partial P'_i=g(h(\partial P_i))$ for each $i=1,\dots,m$. Define $f(P_i):=P'_i$ for each $i$.

Likewise, let $Q_1,\ldots,Q_n$ be the black components of $S^2-M$. By the ``only if'' part of Embedding Extension Theorem there exist disjoint spheres with holes $Q'_1,\ldots,Q'_n$ whose interiors are {\it outside} $g(S^2)$ and such that $\partial Q'_j=g(h(\partial Q_j))$ for each $j=1,\dots,n$. Define $f(Q_j):=Q'_j$ for each $j$.

Image of $f$ is the sphere $\bigsqcup P'_i \sqcup \bigsqcup Q'_j$. Clearly, $f$ and $g$ realize $h$.
\end{proof}

\begin{proof}[Proof of the ``only if'' part of Embedding Extension Theorem]
Consider a properly embedded in $D^3$ sphere with holes $P_i$. Add a ``cap'' (homeomorphic to a disk) in ${\mathbb R^3}-D^3$ to every circle of $\partial P_i$ such that the union of $P_i$ with these caps is a sphere $\hat{P_i}$.
(In the smooth category we may assume that $S^2$ is a round sphere and that bounding circles of $\partial P_i$ are round circles, none of them being an equator. Then for each circle of $\partial P_i$ take the round sphere passing through this circle and the center of $S^2$.
Take parts of such spheres lying in ${\mathbb R^3}-D^3$ as these ``caps''. Analogous, albeit slightly more complicated, construction is possible in the PL category).

Clearly, all same colored connected components of $S^2-M_i=S^2-\partial P_i$ lie on the same side of $\hat{P_i}$. 
 And since $P_i$ and $P_j$ are disjoint, $S^2\cap P_j=M_j$ lie on one side of $\hat{P_i}$, i.e. in the union of same colored components of $S^2-M_i$. 

So $M_j$ lie on one side of $M_i$ by definition. Likewise, $M_i$ lie on one side of $M_j$. Therefore $M_i$ and $M_j$ are unlinked.
\end{proof}

To prove the ``if'' part we require the following claim. Proof of the claim is postponed.
\begin{clm1}
\label{clm1}
Let $P_1,\ldots,P_n$ be properly embedded in $D^3$ pairwise disjoint spheres with holes. Let $M$ be a union of disjoint circles in $S^2=\partial D^3$ such that $M$ and $\partial P_i$ are unlinked for every $i$. Then $M$ lies in one connected component of $D^3-(P_1\sqcup\dots\sqcup P_{n})$.
\end{clm1}

\begin{proof}[Proof of the ``if'' part of Embedding Extension Theorem]
This proof was suggested by A. Novikov. It is simpler than our original proof.

Use induction on number of circles in $M_1\sqcup\dots\sqcup M_m$.

Let $p$ be a circle of $M_1\sqcup\dots\sqcup M_m$ bounding an open disk $D$ in $S^2$ disjoint with $M_1\sqcup\dots\sqcup M_m$ ($p$ corresponds to an edge of $G(S^2,M_1\sqcup\dots\sqcup M_m)$ issuing out of a leaf vertex). We may assume that $p\subset M_1$. Denote by $M'_1$ the union of circles $M_1-p$ (note that $M'_1$ may be empty).

Unions $M'_1, M_2, \ldots, M_m$ are pairwise unlinked. By the inductive hypothesis there are properly embedded in $D^3$ disjoint spheres with holes $P'_1,P_2,\ldots,P_m$ such that $\partial P_i=M_i$ for each $i=2,\dots,m$ and $\partial P'_1=M'_1$. Let $D'$ be a disk obtained from the closure of $D$ by a slight
deformation so that the interior of $D'$ is in the interior of $D^3$ and
$\partial D'=p$. By Claim \ref{clm1} each two points of $M_1=M'_1\sqcup p$ can be connected by a path in the interior of $D^3$ disjoint with $P_2,\dots,P_m$.
So we can connect $D'$ with $P'_1$ by a tube in the interior of $D^3$ disjoint with $P_2,\dots,P_m$. Then we obtain a sphere with holes.
Denote it by $P_1$.
We have $\partial P_1=p\sqcup \partial P'_1=M_1$, $P_1$ is properly embedded in $D^3$ and $P_1$ is disjoint with $P_2,\dots,P_m$.
The inductive step is proved.
\end{proof}

	\begin{figure}[H]
	\centering
	\includegraphics[width=50mm]{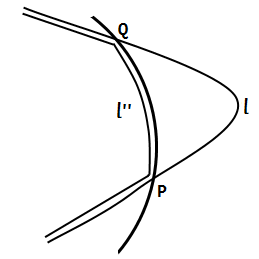}
	\caption{Paths $l$, $l''$.}
	\label{fig:Path}
	\end{figure}
	
\begin{proof}[Proof of Claim \ref{clm1}]
Take any two points $A,B\in M$.
Denote by $l$ a path in $D^3$ connecting $A$ and $B$ such that
$\overline l:=\#(l\cap\bigsqcup\limits_{i=1}^{n} P_i)$ is minimal
(minimal by $l$, objects $A,B,M,D^3,P_1,\dots,P_{n}$ are fixed).
Assume to the contrary that $l$ is not as required, i.e. $\overline l>0$.
Since $M$ is on one side of $\partial P_i$, number $\#(l\cap P_i)$ is even for each $i$.
(If $m=2$, we may even obtain that $\#(l\cap P_1)=0$ and stop here.)
Then $\#(l\cap P_i)\geq 2$ for some $i$.
Denote by $Q$ and $R$ two consecutive points of $l\cap P_i$.
Denote by $Q'$ the point of $l$ slightly before $Q$ and by $R'$ the point of $l$ slightly after $R$ (Fig.~\ref{fig:Path}).
Since $P_i$ is connected, $Q$ and $R$ can be connected by a path in $P_i$.
So $Q'$ and $R'$ can be connected by a path $l'$ very close to $P_i$ but not intersecting $P_i$.
Path $l'$ does not intersect any of $P_1,\ldots,P_n$ because it is very close to $P_i$ and
$P_1,\ldots,P_n$ are pairwise disjoint.
Substitute the part of $l$ between $Q'$ and $S'$ by $l'$.
Denote the obtained path by $l''$.
Then $\overline{l''}=\overline l-2$.
This contradicts to the minimality of $\overline l$.
Thus $l$ is as required.
\end{proof}

\subsection*{Acknowledgements}
I thank A. Novikov for suggestion of a simpler proof of the ``if'' part of Embedding Extension Theorem. I also thank prof. A. Skopenkov for his useful suggestions and comments.

\newpage
{}


\begin{thebibliography}{1}

\bibitem{Lando}
\label{ref:Lando}
I. Arzhantsev, V. Bogachev, A. Garber, A. Zaslavsky, V. Protasov and A. Skopenkov,
\emph{Students' mathematical olympiades at Moscow State University 2010-2011}, Mat. Prosveschenie,
16 (2012), 214-227, in russian, \url{http://www.mccme.ru/free-books/matpros/mpg.pdf}

\bibitem{LKTG}
\label{ref:LKTG}
S. Avvakumov, A. Berdnikov, A. Rukhovich and A. Skopenkov, \emph{How do curved spheres intersect in 3-space,
or two-dimensial meandra}, \url{http://www.turgor.ru/lktg/2012/3/3-1en_si.pdf}

\bibitem{Rukhovich}
\label{ref:Rukhovich}
A. Rukhovich, \emph{On intersection of two embedded spheres in 3-space}, \url{http://arxiv.org/abs/1012.0925}

\bibitem{Hirasa}
\label{ref:Hirasa}
T. Hirasa, \emph{Dissecting the torus by immersions}, Geometriae Dedicata, 145:1 (2010), 33-41

\bibitem{Nowik}
\label{ref:Nowik}
T. Nowik, \emph{Dissecting the 2-sphere by immersions}, Geometriae Dedicata 127, (2007), 37-41, \url{http://arxiv.org/abs/math/0612796}

\end{thebibliography}
\end{document}